\documentclass[noamsfont]{amsproc}

\usepackage{amsmath}
\usepackage{amssymb}
\usepackage{amsthm}
\usepackage{aspmproc}

\newcommand{\ZZ}{\mathbb{Z}}
\newcommand{\Rm}{\mathbb{R}}
\newcommand{\RR}{\mathbb{R}}

\newcommand{\NN}{\mathbb{N}}

\newcommand{\mL}{\mathcal{L}}

\newcommand{\Nm}{\ensuremath{\mathbb{N}}}
\newcommand{\Zm}{\ensuremath{\mathbb{Z}}}

\newcommand{\mM}{\ensuremath{\mathcal{M}}}
\newcommand{\mK}{\ensuremath{\mathcal{K}}}

\newcommand{\mA}{\ensuremath{\mathcal{A}}}

\newtheorem{lem}{Lemma}
\newtheorem{thm}[lem]{Theorem}

\newtheorem{prop}[lem]{Proposition}

\def\qed {\mbox{}\hfill {\small \fbox{}} \\}  
\def\lto{\rightarrow}
\def\lmto{\longmapsto}

\def\leq{\leqslant}
\def\geq{\geqslant}

\title[Weak KAM pairs and Monge-Kantorovich duality]
{Weak KAM Pairs and Monge-Kantorovich Duality}

\author[P.  Bernard and B.  Buffoni]{Patrick Bernard and Boris Buffoni}
\address{
{\rm Patrick Bernard}\\
CEREMADE\\
Universit\'e de Paris Dauphine\\
Pl. du Mar\'echal de Lattre de Tassigny\\
75775 Paris Cedex 16\\
France
\\
\\
{\rm Boris Buffoni}\\
School of Mathematics\\
\'Ecole Polytechnique F\'ed\'erale-Lausanne\\
SB/IACS/ANA Station 8\\
1015 Lausanne\\
Switzerland
}
\begin{document}
\begin{abstract}
The dynamics of globally minimizing orbits of 
Lagrangian systems can be studied using the Barrier function,
as Mather first did, or using the pairs of  weak KAM solutions
introduced by Fathi.
The central observation of the present paper is that
Fathi weak KAM pairs are precisely the admissible pairs for
the Kantorovich  problem dual to the Monge transportation 
problem with the Barrier function as cost.
We exploit this observation to recover
several relations between the Barrier functions and the set of 
weak KAM pairs in an axiomatic and elementary way.
\end{abstract}

\maketitle

\section{Introduction}
Let $M$ be a compact connected manifold and
consider a $C^2$ Lagrangian function
$$
L:TM\times \Rm\lto \Rm
$$
that satisfies the standard hypotheses of the calculus of 
variations,
\begin{equation}\tag{L1}
L(x,v,t+1)=L(x,v,t) 
\ \hbox{ on }\ TM \times \RR,
\end{equation}
\begin{equation}\tag{L2}
\partial^2_{vv}L(x,v,t)>0
\ \hbox{ on }\ TM \times \RR,
\end{equation}
 \begin{equation}\tag{L3}
\lim_{\|v\|\lto \infty} L(x,v,t)/\|v\|=\infty
\ \hbox{ on }M \times \RR.
\end{equation}
 It is standard that, under these assumptions, there exists
 a well-defined time-periodic continuous vectorfield $E(x,v,t)$
 on $TM$ such that the integral curves of $E$ satisfy the
 Euler-Lagrange equations associate to $L$.
 We assume in addition that this vectorfield
generates a complete flow, and denote by $\varphi$
the time-one flow, which is a diffeomorphism
of $TM$.

In this paper we show that the theory developed by Mather 
\cite{Ma:91},  Ma\~n\'e \cite{Mane} and Fathi \cite{Fa:un}
amounts for a large part to the analysis of the function
$A:M\times M\lto \Rm$ defined 
 by the expression
 $$
 A(x,y)=\min_{\gamma}
 \int_0^1 L(\gamma(t),\dot \gamma(t),t)dt,
 $$
 where the minimum is taken on the set of $C^2$ curves
 $\gamma:[0,1]\lto M$ which satisfy 
 $\gamma(0)=x$ and $\gamma(1)=y$.

To emphasize this point of view, we develop an abstract theory
based solely on an arbitrary continuous function $A:M\times M\lto \Rm$,
where $M$ is a connected compact metric space.
We then define $A_1=A$ and
$$A_n(x,y)=\min_{z_1,\ldots,z_{n-1}\in M}A(x,z_1)+A(z_1,z_2)
+\ldots+A(z_{n-1},y)$$
for all integers $n\geq 2$. It turns out that the family
$(A_n)$ is equicontinuous and
our only hypothesis on $A$  is that the family $(A_n)$ is uniformly
bounded (this can be achieved by adding some constant to $A$). 
It then follows that the expression
$$c(x,y)=\liminf_{n\rightarrow \infty}A_n(x,y)$$ 
defines a continuous function $c:M\times M\lto \RR$.

We call $(\phi_0,\phi_1)$ an {\em admissible Kantorovich pair}
for $c$ if
$$
\forall y\in M\ \ \phi_1(y)=\min_{x\in M} \phi_0(x)+c(x,y)
$$
 and  
$$
\forall x\in M\ \ \phi_0(x)=\max _{y\in M} \phi_1(y)-c(x,y).
$$
The first main result (Theorem \ref{pairs}) states
that $(\phi_0,\phi_1)$ is an admissible Kantorovich pair
for $c$ if and only if
\begin{itemize}
\item $\phi_0(x)=\max_{y\in M}\phi_0(y)-A(x,y)$ for all $x\in M$,
\item $\phi_1(x)=\min_{y\in M}\phi_1(y)+A(y,x)$ for all $x\in M$,
\item and $\phi_0(x)=\phi_1(x)$ whenever $c(x,x)=0$.
\end{itemize}

The second main result 
(Theorem \ref{Mather})
concerns the minimization problem
$$\min_{\eta} \int_{M\times M} A(x,y)d\eta(x,y),$$
where the minimum is taken on the set of Borel
probability measures $\eta$ on $M\times M$ with equal marginal measures, 
that is,
$\pi_{0\sharp}(\eta)=\pi_{1\sharp}(\eta)$ with $\pi_0$ and $\pi_1$
denoting the canonical projections on $M$.
Among all admissible measures, the minimizing ones
are shown to be exactly those supported on the set
$$D=\{(x,y)\in M\times M\,|\, A(x,y)+c(y,x)=0\}.$$

This is also restated in the following  way in Theorem \ref{Mather2}. 
Let $X=M^{\ZZ}$
be endowed with the product topology and denote by $\mM_T(X)$
the set of Borelian probability measures on $X$
which are invariant by translation. Consider the
minimization problem
$$\min_{\nu\in \mM_T(X)}\int_X A(x_0,x_1)d\nu(x),$$
where a generic $x\in X$ is written $x=(\ldots,x_{-1},x_0,x_1,\ldots)$.
Then we show with the help of the Ergodic Decomposition
Theorem that $\nu$ in $\mM_T(X)$ is minimizing exactly when
the push-forward of $\nu$ by the projection 
$x\rightarrow (x_0,x_1)$
is concentrated on $D$.

The paper ends with the interpretation of these abstract theorems
in the setting of the Aubry-Mather theory, recovering in this way
some key results of
\cite{Mane,Ma:91,Fa:un}.

\section{Monge-Kantorovich theory}\label{MKT}
We present some standard facts of Monge-Kantorovich 
theory, first in the general case, and then 
when the cost satisfies some given assumptions.
\subsection{Generalities}
We recall the basics of Monge-Kantorovich duality.
The proofs are available in many texts on the subjects, for
example \cite{Am:00,RaRu:98,Vi:03}.
We assume that $M$ and $N$ are compact metric spaces, and that 
$c(x,y)$ is a continuous cost function on 
$M\times N$.
Given  Borel  probability measures $\mu_0$ on $M$ 
and $\mu_1$ on $N$,
a transport plan between $\mu_0$ and 
$\mu_1$ is a measure on $M\times N$
which satisfies 
$$\pi_{0\sharp}(\eta)=\mu_0 \text{ and }
\pi_{1\sharp}(\eta)=\mu_1,
$$
where $\pi_0:M\times N\lto M $ is the projection
on the first factor, and $\pi_1:M\times N\lto N$ is the projection on
the second factor. 
We denote by $\mK(\mu_0,\mu_1)$, after Kantorovich,
the set of transport plans.
Kantorovich proved the existence of 
a minimum in the expression 
\begin{equation}\label{Monge}
C(\mu_0,\mu_1)=
\min_{\eta\in \mK(\mu_0,\mu_1)} \int _{M\times N} c d\eta
\end{equation}
for each  pair $(\mu_0, \mu_1)$ of probability measures.
The plans which realize this minimum are called
optimal transfer plans.
Let $\phi_0$ be a real function on $M$ and $\phi_1$ 
a real function on $N$.
The  pair $(\phi_0,\phi_1)$ 
is called an admissible Kantorovich pair if it satisfies the 
relations
$$
\phi_1(y)=\min_{x\in M} \phi_0(x)+c(x,y)
\text{ and } 
\phi_0(x)=\max _{y\in N} \phi_1(y)-c(x,y)
$$
for all point  $x\in M$ and $y\in N$.
Another discovery of Kantorovich is that 
\begin{equation}\label{Kanto}
C(\mu_0,\mu_1)
=\max _{\phi_0,\phi_1}
\Big(
\int_N \phi_1d\mu_1-\int_M \phi_0d\mu_0
\Big)
\end{equation}
where the maximum is taken on the non-empty  set of 
 admissible Kantorovich pairs $(\phi_0, \phi_1)$.
This maximization problem is called the dual Kantorovich problem,
the admissible pairs  which reach this maximum
are called optimal Kantorovich pairs.
The direct problem (\ref{Monge}) and dual problem (\ref{Kanto})
are related as follows. 
\begin{prop}\label{fathi}
If $\eta$ is an optimal transfer plan, and if $(\phi_0,\phi_1)$
is a Kantorovich  optimal pair, then the support of $\eta$ 
is contained in the set 
$$
\{(x,y)\in M\times N\text{ such that }
\phi_1(y)-\phi_0(x)=c(x,y)\},
$$
which is a closed subset of $M\times N$ 
because $\phi_0$ and $\phi_1$ are continuous.
\end{prop}
Let us remark that the knowledge of the set of 
Kantorovich admissible pairs
is equivalent to the knowledge of the cost function $c$.
\begin{lem}\label{var char of c}
We have 
$$c(x,y)=\max _{(\phi_0,\phi_1)} \phi_1(y)-\phi_0(x)
$$
where the maximum is taken on the set of Kantorovich
admissible pairs.
\end{lem}
\proof
This Lemma is elementary and can be proved by easy manipulation
of inequalities, see \cite{BeBu:first}.
 However, we present a short  proof based on the
non-elementary Monge-Kantorovich duality.
Let us
fix  points $x\in M$ and $y\in N$, and let $\mu_0$ be the Dirac measure
at $x$ and $\mu_1$ be the Dirac measure at $y$.
There exists one and only one transport plan between
$\mu_0$ and $\mu_1$, it is the Dirac measure at $(x,y)$.
As a consequence, we have $c(x,y)=C(\mu_0,\mu_1)$.
Hence the equality above is precisely the conclusion
of Kantorovich duality for the transportation problem between
$\mu_0$ and $\mu_1$.
\qed
\begin{prop}\label{converse}
Let $(\phi_0,\phi_1)$ be an admissible pair,
and let $\mu_0$ be a probability measure on $M$.
Then there exists a probability  measure $\mu_1$ on $N$ such that 
the pair $(\phi_0,\phi_1)$ is optimal for the transportation
problem of the measure $\mu_0$ onto the measure $\mu_1$
\end{prop}
\proof
If $\mu_0$ is the Dirac at $x$, then take a point $y$
such that $\phi_1(y)=\phi_0(x)+c(x,y)$, and observe that
the conclusion obviously holds if $\mu_1$
is the Dirac at $y$.
The set of measures $\mu_0$ for which the conclusion
holds (given $\phi_0,\phi_1$) is clearly convex and closed 
(with respect to the weak topology), 
it contains the Dirac measures,
hence it is the whole set of probability measures.
\qed

\subsection{Distance-like costs}
Kantorovich stated his duality theorem
first in the case where $M=N$ and the cost is a distance.
Then, the dual problem  takes a simpler
 form that we now describe.
In fact, it is not necessary to assume that the cost is a distance.
It is sufficient to assume that, for all $x$, $y$ and $z$ in $M$,
we have
\begin{equation}\tag{C1}
c(x,z)
\leq c(x,y)+c(y,z),
\end{equation}
\begin{equation}\tag{C2}
c(x,x)= 0.
\end{equation}
A function $\phi:M\lto \Rm$ is called $c$-Lipschitz
if it satisfies the inequality
$$
\phi(y)-\phi(x)\leq c(x,y)
$$
for all $x$ and $y$ in $M$.
Note that, in the above and in what follows, we assume 
that $M=N$ is a compact
and connected metric space, and that $c:M\times M\lto \Rm$
is a continuous cost function.
\begin{thm}
Assume that the cost $c\in C(M^2,\Rm)$ satisfies the assumptions 
$(C1)$ and $(C2)$. Then for each pair $\mu_0$, $\mu_1$
of probability measures on $M$, we have
$$
C(\mu_0,\mu_1)=
\max _{\phi} \int_M \phi d(\mu_1-\mu_0)
$$
where the maximum is taken on the set of $c$-Lipschitz 
functions $\phi$.
\end{thm}
This is a 
well-known
direct rewriting of Kantorovich duality in view
of the following description of admissible pairs.
\begin{lem}\label{dist}
If the cost satisfies $(C1)$ and $(C2)$, then
the Kantorovich admissible pairs are precisely
the pairs of the form $(\phi,\phi)$, with $\phi$
$c$-Lipschitz.
\end{lem}
\proof
If $\phi$ is a $c$-Lipschitz function,
then $(\phi,\phi)$ is an admissible pair.
Indeed, let us prove for example that
$\phi(x)=\min_y \phi(y)+c(y,x)$.
On the one hand, we have $\phi(x)\leq \phi(y)+c(y,x)$
because $\phi$ is $c$-Lipschitz, hence 
$\phi(x)\leq \min_y \phi(y)+c(y,x)$.
On the other hand, 
$\phi(x)=\phi(x)+c(x,x)\geq \min_y \phi(y)+c(y,x)$.
One can prove similarly that
$\phi(x)=\max_y \phi(y)-c(x,y)$.
It follows that $(\phi,\phi)$ is an admissible pair.
Conversely, if 
$(\phi_0,\phi_1)$ is an admissible pair,
then $\phi_0=\phi_1$ is a $c$-Lipschitz function.
This is a special case of Lemma \ref{clip} below.
\qed

Let us now study costs which satisfy $(C1)$ but not necessarily
$(C2)$.
It is then useful to define the set
$$
\mA:=\{x\in M, c(x,x)=0\}\subset M.
$$
Note that the restriction of the cost $c$ to $\mA\times \mA$
obviously satisfies $(C1)$ and $(C2)$.
In this more general case, we have:
\begin{lem}\label{clip}
Let $c\in C(M^2,\Rm)$ satisfy $(C1)$.
Let $(\phi_0,\phi_1)$
be an admissible pair.
Then the functions $\phi_0$ and $\phi_1$ 
are $c$-Lipschitz.
In addition, we have $\phi_0\leq \phi_1$ with equality on $\mA$.
\end{lem}

\proof
Let us first prove that the function 
$\phi_1$
is $c$-Lipschitz.
Given $x\in M$, there exists $y$ such that 
$
\phi_1(x)=\phi_0(y)+c(y,x)
$,
and then, for each $z$,
$$
\phi_1(x)=\phi_0(y)+c(y,x)\geq \phi_1(z)-c(y,z)+c(y,x)
\geq \phi_1(z)-c(x,z).
$$
One can prove similarly that $\phi_0$ is $c$-Lipschitz.

We then have 
$$
\phi_0(x)=\max _y \phi_1(y)-c(x,y)
\leq \max _y \phi_1(x)=\phi_1(x).
$$
because $\phi_1$ is $c$-Lipschitz.
If $x\in \mA$, we have, in addition,
$$
\phi_0(x)=\max _y \phi_1(y)-c(x,y)
\geq  \phi_1(x)-c(x,x)=\phi_1(x).
$$
\qed
We now introduce another hypothesis which
is certainly less natural than $(C1)$ and
$(C2)$, but is useful for the applications we have in mind.
We assume that
\begin{equation}\tag{C3}
\mA\neq\emptyset \ \ \hbox{ and }\ \ 
c(x,y)=\min_{a\in \mA} c(x,a)+c(a,y)
\end{equation}
for  each $x$ and $y$ in $M$. 
Note that, under  $(C1)$, $(C3)$ is implied by $(C2)$.
The hypothesis $(C3)$ implies that 
 each optimal transport can be factored
through the set $\mA$.

\begin{lem}\label{factor}
If the cost satisfies $(C1)$ and $(C3)$, then
for each pair $(\mu_0,\mu_1)$ of probability measures,
there exists a probability measure $\mu$ 
 supported on 
$\mA$ and such that
$$
C(\mu_0,\mu_1)=C(\mu_0,\mu)+
C(\mu,\mu_1)
$$
\end{lem}
\proof
First note that
$C(\mu_0,\mu_1)\leq C(\mu_0,\mu)+
C(\mu,\mu_1)$ is true for all Borelian probability measures $\mu$
on $M$. This can be seen as follows. Let $\eta_0$ and $\eta_1$
be optimal transport plans for $(\mu_0,\mu)$ and $(\mu,\mu_1)$ respectively.
Disintegrate $\eta_0$ with respect to $\pi_1$ and $\eta_1$ with
respect to $\pi_0$: $\eta_0=\int_M\eta_{0z}d\mu(z)$ and
$\eta_1=\int_M\eta_{1z}d\mu(z)$ 
(see e.g. Theorem 5.3.1 in \cite{AmGiSa} for the disintegration theorem;
here $\eta_{0z}$ and $\eta_{1z}$ are seen as probability measures on $M$).
Following Section 5.3 in \cite{AmGiSa}, 
define the probability measure $\eta$ on $M^2$ by
$$\eta(A\times B)=\int_M \eta_{0z}(A)\eta_{1z}(B)\, d\mu(z)$$
for all Borelian subsets $A,B\subset M$. Then $\eta\in\mK(\mu_0,\mu_1)$
and 
\begin{multline*}
\int_{M^2} c\, d\eta
=\int_{M^3} c(x,y)d\eta_{0z}(x)d\eta_{1z}(y) d\mu(z)
\\ \leq\int_{M^3} \{c(x,z)+c(z,y)\}d\eta_{0z}(x)d\eta_{1z}(y) d\mu(z)
= \int_{M^2} c\, d\eta_0+\int_{M^2} c\, d\eta_1.
\end{multline*}

Let us now prove the reverse inequality  when $\mu_0$ and $\mu_1$
are Dirac measures supported in $x$ and $y$.
In this case, one can take for $\mu$ the Dirac measure
supported at $a$, where $a$ is any point such that
$c(x,y)=c(x,a)+c(a,y)$.
The general case is then deduced once again
using the fact that, 
on $M^2$, the set of probability measures
is the closed convex envelop of the set of Dirac measures, so that
we can approximate any optimal transfer plan in $\mK(\mu_0,\mu_1)$
by Dirac measures.
\qed

\begin{prop}\label{admissiblepairs}
If the cost $c\in C(M^2,\Rm)$ satisfies $(C1)$ and $(C3)$, then for each 
admissible pair $(\phi_0,\phi_1)$, 
there exists a 
function $\phi:\mA\lto \Rm$, which is $c$-Lipschitz, and such that
$$
\phi_1(a)=\phi_0(a)=\phi(a)
$$
for all $a\in \mA$,
\begin{equation}\label{phi1}
\phi_1(x)=\min_{a\in \mA} \phi(a)+c(a,x)
\end{equation}
for all $x\in M$
and 
\begin{equation}\label{phi0}
\phi_0(x)=\max_{a\in \mA} \phi(a)-c(x,a).
\end{equation}
Conversely, given any $c$-Lipschitz function $\phi$ on $\mA$,
the functions $\phi_0$ and $\phi_1$ defined by 
(\ref{phi0}) and (\ref{phi1}) form  an admissible pair.
In other words, there is a bijection between the set
of admissible pairs and the set of $c$-Lipschitz functions on $\mA$.

\end{prop}
\proof
The fact that $\phi_0$ and $\phi_1$ are
$c$-Lipschitz and that, on $\mA$,  
$\phi_0=\phi_1:=\phi$ results from Lemma \ref{clip}.
Let us prove (\ref{phi0}), the proof of (\ref{phi1})
being similar:
$$
\phi_0(x)=\max _y \phi_1(y)-c(x,y)
\stackrel{(C3)}{=}
\max _{y\in M,a\in \mA} \phi_1(y)-c(x,a)-c(a,y)
$$
$$
=
\max _{a\in \mA} \phi_0(a)-c(x,a)
=
\max _{a\in \mA} \phi(a)-c(x,a).
$$

Conversely, let $\phi$ be a $c$-Lipschitz function
on $\mA$, and let $\phi_0$ and $\phi_1$
be defined by (\ref{phi0}) and (\ref{phi1}).
The reader will easily check that $\phi_0$ 
and $\phi_1$ are $c$-Lipschitz,
and that $\phi_1\leq \phi\leq \phi_0$ on $\mA$.
We now prove that $\phi_0\leq \phi_1$ (and then that there is 
equality on $\mA$):
$$
\phi_0(x)-\phi_1(x)
=
\max_{a,b\in \mA} \phi(a)-c(x,a)-\phi(b)-c(b,x)
$$
$$
\leq \max_{a,b\in \mA} \phi(a)-\phi(b)-c(b,a)\leq 0
$$
because $\phi$ is $c$-Lipschitz on $\mA$.
In order to check that the pair $(\phi_0,\phi_1)$
is an admissible pair,
we shall prove that 
$$
\phi_0(x)=\max _y \phi_1(y)-c(x,y)
$$
and 
leave the other half to the reader.
For each $x$  in $M$, we have 
$$
\phi_0(x)=\max_{a\in \mA} \phi(a)-c(x,a)
=\max_{a\in \mA} \phi_1(a)-c(x,a)\leq
\max_{y\in M} \phi_1(y)-c(x,y).
$$
In order to obtain the other inequality,
let us prove that 
$$
\phi_1(y)-\phi_0(x)\leq c(x,y)
$$
for all $x$ an $y$ in $M$.
Indeed, we have
$$
\phi_1(y)-\phi_0(x)=
\min_{a,b\in \mA} \phi(a)+c(a,y)-\phi(b)+c(x,b)
$$
$$
\leq \min_{a,b\in \mA} c(b,a)+c(a,y)+c(x,b)
= \min_{a\in \mA} c(x,a)+c(a,y)=c(x,y)
$$
by $(C3)$.
\qed

Since $\phi_1$ and $\phi_2$ are $c$-Lipschitz (Lemma \ref{clip}),
equations \eqref{phi1} and \eqref{phi0}  imply
$$
\phi_1(x)=\min_{y\in M} \phi_1(y)+c(y,x)
\ \hbox{ and }  \ 
\phi_0(x)=\max_{y\in M} \phi_0(y)-c(x,y)
$$
for all $x\in M$.

\section{Abstract Mather-Fathi Theory}

In this section, we consider a continuous
function $A(x,y):M\times M\lto \Rm$.
Recall that $M$ is a compact connected metric space.
We shall build several functions out of $A$.
First, we define the sequence of functions $A_n(x,y)$
by  setting $A_1=A$ and
\begin{multline*}
A_n(x,y)=
\min _{z\in M}
A(x,z)+A_{n-1}(z,y)
\\=\min _{z_1,\ldots,z_{n-1}\in M} A(x,z_1)+A(z_1,z_2)+\ldots + A(z_{n-1},y).
\end{multline*}

 \begin{lem}\label{asymptotique}
 The functions $A_n$ are equicontinuous.
 In addition, there exists a real number $l$ 
 and a positive constant $C$ such that
 $$
 |A_n(x,y) -ln|\leq C 
 $$ 
 for all $n \in \Nm$
 and all $x$ and $y$ in $M$.
 \end{lem}
 
 \proof
The function $A$ is continuous, hence uniformly continuous,
hence there exists a modulus
of continuity $\delta:[0,\infty)\lto [0,\infty)$
such that $\lim_{\epsilon\lto 0} \delta(\epsilon)=\delta(0)=0$
and such that
$$
|A(x,y)-A(X,Y)|\leq \delta(d(x,X))+\delta(d(y,Y))
$$
for all $x,y,X,Y$ in $M$.
Clearly, for all $n\geq 2$ and all $z_1,\ldots,z_{n-1}\in M$, the function
$(x,y)\rightarrow  A(x,z_1)+A(z_1,z_2)+\ldots + A(z_{n-1},y)$ is uniformly continuous,
with the same modulus of continuity as $A$.
Hence $A_n$ is uniformly continuous with the same modulus of continuity as $A$ because it is the infimum
of functions having all the same modulus of continuity.

Let us define the sequences $M_n:=\max_{(x,y)\in M^2}A_n(x,y)$
and $m_n:=\min_{(x,y)\in M^2}A_n(x,y)$.
It is clear that the sequence $M_n$
is subadditive, \textit{i.~e.}
that $M_{n+k}\leq M_n+M_k$ for all $n$ and $k$ in $\Nm$.
In order to check this claim, we take $x$ and $y$ in $M$  such that
$A_{n+k}(x,y)=M_{n+k}$.
Then there exists a point $z$ in $M$ such that
$$
M_{n+k}=A_{n+k}(x,y)=A_n(x,z)+A_k(z,y)\leq M_n+M_k.
$$
Similarly,
the sequence $m_n$ is super-additive, \textit{i. e.}
 $m_{n+k}\geq m_n+m_k$.
On the other hand, on view of the equicontinuity of 
$A_n$, there exists a constant $C$ such that 
$M_n-m_n\leq C$.
Applying a standard result on subadditive sequences (see e.g.
Lemma 1.18 in \cite{Bo}), we obtain
that $M_n/n$ converges to its infimum $M$,
and that $m_n/n$ converges to its supremum $m$.
Then for each $x$ and $y$,
 $$
 nM-C\leq M_n-C\leq m_n\leq A_n(x,y)
 \leq M_n\leq  m_n+C\leq nm+C
 $$
which implies that $M=m$ and proves the Lemma.
\qed

We make, on the function $A$, the hypothesis
\begin{equation}\tag{A1}
l=0.
\end{equation}
Note that this hypothesis implies that 
$A(x,x)\geq 0$ for all $x$, and more generally
that $A_n(x,x)\geq 0$ for all $x$.
Then, we can  define a cost function $c$ by the expression

\begin{equation}\label{defc}
c(x,y)=\liminf_{n\rightarrow \infty} A_n(x,y).
\end{equation}

In view of Lemma \ref{asymptotique},
the function $c$ takes finite values and is continuous.
We have $c(x,x)\geq 0$ and,
by Lemma \ref{lem: check C3} below, $c(x,y)+c(y,x)\geq c(x,x)\geq 0$
for all $x$ and $y$ in $M$.

\begin{lem}\label{propdec}
For each $n\in\Nm$, 
we have
$$
c(x,y)=\min_{z\in M} c(x,z)+A_n(z,y)=
\min_{z\in M} A_n(x,z)+c(z,y).
$$
\end{lem}

\proof
Let us fix $n$.
Passing at the liminf ($m\lto \infty$) in the inequality
$$
A_{m+n}(x,y)\leq A_m(x,z)+A_n(z,y),
$$
we obtain
$$
c(x,y)\leq c(x,z)+ A_n(z,y).
$$
For the opposite inequality, let us notice that,
for each $m$, there exists a point $z_m$
in $M$ such that 
$$
A_{m+n}(x,y)
=A_m(x,z_m)+A_n(z_m,y).
$$
Let us consider an increasing  sequence of integers $m_k$
such that the subsequence $z_{m_k}$ has a limit $z$ and 
$\lim_{k\rightarrow \infty} A_{m_k+n}(x,y)=c(x,y)$.
At the liminf, we get, taking advantage of the equicontinuity of the 
functions $A_n$,
$$
c(x,y)\geq c(x,z)+A_n(z,y).
$$
This proves that 
$$
c(x,y)=\min_z c(x,z)+A_n(z,y).
$$
The proof of the second equality of the statement is similar.
\qed

\begin{lem}\label{lem: check C3}
The cost function $c$ satisfies $(C1)$
and $(C3)$.
\end{lem}
\proof
 The triangle inequality is easily deduced from
 Lemma \ref{propdec}.
 Let us now prove $(C3)$.
 We first prove that, given $x$ and $y$ in $M$,
 there exists a point $z$ in $M$ such that 
 $c(x,y)=c(x,z)+c(z,y)$.
 Indeed, for each $n$ in $\Nm$,
 there exists  
  a point $z_n$ such that
  $c(x,y)=c(x,z_n)+A_n(z_n,y)$.
 Considering an increasing sequence of integers $n_k$
 such that the subsequence $z_{n_k}$ has a limit $z$,
 we obtain at the liminf along this subsequence that
 $c(x,y)\geq c(x,z)+c(z,y)$
 which is then an equality.
 
 By recurrence, there exists a sequence $Z_n\in M$
 such that, for each $k\in \Nm$,
 we have 
 $$
 c(x,y)=c(x,Z_1)+c(Z_1,Z_2)+\ldots+ c(Z_{k-1},Z_k)
 +c(Z_k,y).
 $$  
Note that $\sum_{i=\ell}^{m}c(Z_i,Z_{i+1})
=c(Z_{\ell},Z_{m+1})$ if $0\leq \ell <m\leq k$, where $Z_0=x$
and $Z_{k+1}=y$.

Let $Z$ be an accumulation point of the sequence $Z_n$.
For each $\epsilon>0$, we can suppose, by taking a subsequence 
in $Z_n$, that all the points $Z_n$ belong to the ball
of radius $\epsilon$ centered at $Z$.
We conclude that, for each $k\in \Nm$,
$$
c(x,y)\geq c(x,Z)+(k-1)c(Z,Z)+c(Z,y)-2(k+1)\delta(\epsilon).
$$
This is possible only if $c(Z,Z)\leq 2\delta(\epsilon)$,
and since this should hold for all $\epsilon$
we conclude that $c(Z,Z)\leq 0$,
hence $c(Z,Z)=0$.
We have proved the existence of a point $Z\in \mA$
such that $c(x,y)=c(x,Z)+c(Z,y)$.
\qed

Let us define, the two operators $T^{\pm}$  on the
space $C(M,\Rm)$ of continuous functions on $M$ by the
expressions
 $$
T^-u(x)=
\min_{y\in M}
u(y)+ A(y,x)
$$
and
$$
T^+u(x)=
\max_{y\in M}
u(y)- A(x,y).
$$

We have the following relation between the fixed points of these
operators and the admissible pairs of the Kantorovich dual problem
with cost $c$.
Recall the definition 
$
\mA:=\{x\in M, c(x,x)=0\}\subset M.
$

\begin{thm}\label{pairs}
Let $A$ be a function satisfying $(A1)$, and let $c$
be the cost defined by (\ref{defc}).
The pair $(\phi_0,\phi_1)$ of functions on $M$ 
is a Kantorovich admissible pair (for $c$) if and only if
\begin{itemize}
\item the function $\phi_0$ is a fixed point of $T^+$,
\item the function $\phi_1$ is a fixed point of $T^-$,
\item $\phi_0=\phi_1$ on $\mA$.
\end{itemize}
Finally, for each fixed point $\phi_1$ of $T^-$,
there exists one and only one function $\phi_0$
such that $(\phi_0,\phi_1)$ is an admissible pair.
\end{thm}
\proof
Let $(\phi_0,\phi_1)$
be an admissible pair.
Then we have the expression 
$$
\phi_1(y)=
\min_{x\in M} \phi_0(x)+c(x,y).
$$
We obtain that
$$
T^-\phi_1(z)=\min_{x,y\in M} \phi_0(x)+c(x,y)+ A(y,z)
$$
$$
=\min_{x\in M} \phi_0(x)+c(x,z)=\phi_1(z).
$$
We prove in the same way that 
the function $\phi_0$ is a fixed point
of $T^+$.
 Lemma \ref{clip} implies that $\phi_0=\phi_1$ on $\mA$.

Conversely, let $(\phi_0,\phi_1)$ satisfy the three conditions
of the statement. We first observe that the functions $\phi_0$ and
$\phi_1$ are $c$-Lipschitz. 
Indeed, we have, for each $n$, 
$$
\phi_i(y)-\phi_i(x)\leq A_n(x,y).
$$
When $n=1$, this is a direct consequence of the fact that $\phi_i$
is a fixed point of $T^{\pm}$, and the general case is proved by
induction.
We get
$$\phi_i(y)-\phi_i(x)\leq\liminf_{n\rightarrow \infty}A_n(x,y)= c(x,y).$$
The function $\phi_1$ being a fixed point of $T^-$, for each 
$n\in \Nm$, there exists a point $y_n$ in $M$
such that
$\phi_1(x)=\phi_1(y_n)+A_n(y_n,x)$.
Indeed, we can find successively $y_1,y_2,\ldots$ such that
\begin{multline*}
\phi_1(x)=\phi_1(y_1)+A(y_1,x)
=\phi_1(y_2)+A(y_2,y_1)+A(y_1,x)
\\=\ldots=\phi_1(y_n)+A(y_n,y_{n-1})+\ldots+A(y_1,x).
\end{multline*}
By definition of $A_n$, we get $\phi_1(x)\geq \phi_1(y_n)+A_n(y_n,x)$.
The reverse inequality has just been proved above.

Let $n_k$ be a subsequence
 such that $y_{n_k}$ has a limit $y$.
 At the limit, we obtain the inequality
 $$
 \phi_1(x)\geq \phi_1(y)+c(y,x),
 $$
which is then an equality.
We have proved that 
$$
\phi_1(x)=
\min_{y\in M} \phi_1(y)+c(y,x).
$$
Let us call $\phi$ the common value of $\phi_0$ and $\phi_1$ on $\mA$.
In view of $(C3)$, we have
$$
\phi_1(x)=\min_{y\in M, a\in \mA}
\phi_1(y)+c(y,a)+c(a,x)
$$
$$
=
\min_{a\in \mA} \phi_1(a)+c(a,x)
=\min_{a\in \mA} \phi(a)+c(a,x).
$$
One can prove in a similar way that
$$
\phi_0(x)=\max _{a\in \mA} \phi_0(a)-c(x,a)
=\max _{a\in \mA} \phi(a)-c(x,a).
$$
We conclude that $(\phi_0,\phi_1)$ is an admissible pair
by Proposition \ref{admissiblepairs}.
This also proves the uniqueness claim.

In order to prove the last part of the statement,
let us consider a fixed point $\phi_1$ of $T^-$.
Let us define the function $\phi_0$
by
$$
\phi_0(x)=\max _{a\in \mA} \phi_1(a)-c(x,a).
$$
Since the function $\phi_1$ is $c$-Lipschitz (as seen above),
we have $\phi_0\leq \phi_1$.
On the other hand, it is clear that $\phi_1\leq \phi_0$
on $\mA$.
As a consequence, we have $\phi_0=\phi_1$ on $\mA$. 
By Lemma \ref{propdec},
we have  for   all $z\in M$ that
\begin{multline*}
\max_{x\in M}\phi_0(x)-A(z,x)
=\max_{x\in M,a\in \mA}\phi_1(a)-c(x,a)-A(z,x)
\\=\max_{a\in \mA} \phi_1(a)-c(z,a)
=\phi_0(z).
\end{multline*}
Hence the function $\phi_0$
is a fixed point of $T^+$
and, as a consequence, the pair $(\phi_0,\phi_1)$
is an admissible pair.
\qed


\section{Dynamics}\label{dyn}
Let us define the subset 
$$
D:=\{(x,y)\in M\times M
\text{ s. t. } A(x,y)+c(y,x)=0\}\subset \mA\times\mA
$$
(see Lemma \ref{propdec}).
We shall explain in two different ways 
that the Borel probability measures $\eta$
on $M\times M$  which are supported on $D$
and satisfy 
$
\pi_{0\sharp}(\eta)=
\pi_{1\sharp}(\eta)
$
can be seen in a natural way as 
the analog of Mather minimizing measures in our setting.\

\subsection{Construction via Kantorovich pairs}
We first expose a construction based on Kantorovich
pairs.

\begin{thm}\label{Mather}
Under the assumption $(A1)$,
we have 
$$
\min _{\eta}
\int_{M\times M}A(x,y)d\eta(x,y)
=0,
$$
where the minimum is taken on the set of Borel probability
measures 
$\eta$ on $M\times M$ such that 
$
\pi_{0\sharp}(\eta)=
\pi_{1\sharp}(\eta).
$
The minimizing measures  are those which are supported on $D$.
\end{thm}

\proof
Let us first prove that there exists a measure 
$\eta$ on $M\times M$  which is supported on $D$ and such that 
$
\pi_{0\sharp}(\eta)=
\pi_{1\sharp}(\eta).
$
By Lemma \ref{propdec}, for each $x_0\in \mA$,
there exists a point $x_1$ in $\mA$
such that $(x_0,x_1)\in D$.
Hence there exists a sequence 
$x_0, x_1, x_2, \ldots x_n, \ldots$
of points of $\mA$ such that 
$(x_n,x_{n+1})\in D$ for each $n$.
Let us now consider the sequence
$$
\eta_n
=
\frac{\delta_{(x_0,x_1)}+
\delta_{(x_1,x_2)}+\cdots+
\delta_{(x_{n-1},x_n)}}
{n}
$$
of probability measures on $\mA\times \mA$.
Every accumulation point (for the weak topology)
of the sequence $\eta_n$ satisfies the desired property.
Since the set of probability measures on $M\times M$
is compact for the weak topology, such accumulation points exist.

Consider a measure
 $\eta$ on $M\times M$  which is supported on $D$ and such that 
$
\pi_{0\sharp}(\eta)=
\pi_{1\sharp}(\eta).
$
We have 
$$
\int A(x,y)d\eta(x,y)
=
\int -c(y,x)d\eta(y,x)
\leq\int \phi(y)-\phi(x) d\eta(x,y)
=0,
$$
where $\phi$ is any $c$-Lipschitz function.

On the other hand,
let
$\eta$ 
be a probability measure on $M\times M$ such that 
$
\pi_{0\sharp}(\eta)=
\pi_{1\sharp}(\eta).
$
Consider a function $\phi$ which is $A$-Lipschitz.
Such functions exist, for example, take 
$z_2\rightarrow c(z_1,z_2)$ for any $z_1\in M$ (see Lemma~\ref{propdec})
or fixed points of $T^-$ or $T^+$.
We have 
\begin{equation}\label{in}
0= \int \phi(y)-\phi(x) d\eta(x,y)
\leq \int A(x,y)d\eta(x,y).
\end{equation}
We have proved that the minimum in the statement is indeed
zero, and that the measures supported on $D$
are minimizing.
There remains to prove that every minimizing measure is 
supported on $D$.

It is clear that a measure $\eta$ is minimizing if and only if,
for each $A$-Lipschitz function $\phi$, there is equality
in (\ref{in}), which means that the measure $\eta$
is supported on the set 
\begin{multline*}
D_1=\{ (x,y)\in M^2\,|\,
\\ \phi(y)-\phi(x)=A(x,y)
\hbox{ for all $A$-Lipschitz functions $\phi$}\}.
\end{multline*}
Let $D_{\infty}$ be the set of pairs 
$(x_0,x_1)$ such that there exists a sequence
$x_i$, $i\in \Zm$
satisfying $(x_i,x_{i+1})\in D_1$ for all $i\in \Zm$
(and of course with the given points $x_0$ and $x_1$).

We claim  that  $D_{\infty}\subset D$.
In order to prove this claim, let 
$\phi$ be $A$-Lipschitz. Observe that $\phi$ is $A_n$-Lipschitz
for all $n\in \NN$ and $c$-Lipschitz.
If $(x_0,x_1)$ is a point in $D_{\infty} $,
then there exists a sequence $x_i$, $i\in \Zm$
such that 
$$\phi(x_j)-\phi(x_i)= A_{j-i}(x_i,x_j)$$
for each $i< j$ in $\Zm$.
If $\alpha$ is an accumulation point of the sequence 
$x_i$ at $-\infty$, 
we get  the equality
$$\phi(x_j)-\phi(\alpha)= c(\alpha,x_j)$$
for each $j\in \ZZ$ and then, in the same way,
$c(\alpha,\alpha)=\phi(\alpha)-\phi(\alpha)=0$, hence 
$\alpha\in \mA$. Let $(\phi_0,\phi_1)$ be a Kantorovich pair
for $c$, so that both $\phi_0$ and $\phi_1$ are $A$-Lipschitz
(see Theorem \ref{pairs}).
We get
$\phi_1(\alpha)=\phi_0(\alpha)$ (because $\alpha\in \mA$, see Theorem
\ref{pairs})
hence $\phi_1(x_j)=\phi_0(x_j)$.
Since this holds for all Kantorovich pairs, we get that
$x_j\in \mA$ (see Lemma \ref{var char of c}).
In other words, we have proved that 
$D_{\infty} \subset \mA\times \mA$.
Now let $(x_0,x_1)$ be a point of $D_{\infty}$.
We have
$x_1\in \mA$, and, since the function
$c(x_1,.)$ is $A$-Lipschitz,
we have the equality
$
c(x_1,x_1)-c(x_1,x_0)=A(x_0,x_1).
$
Recalling that $c(x_1,x_1)=0$, we get 
$c(x_1,x_0)+A(x_0,x_1)=0$,
hence $(x_0,x_1)\in D.$
The proof of the Theorem then follows from the next Lemma.
\qed

\begin{lem}
If $\eta$ is a probability measure on $M\times M$
which is supported on $D_1$ and such that 
$
\pi_{0\sharp}(\eta)=
\pi_{1\sharp}(\eta),
$
then $\eta$ is concentrated on $D_{\infty}$.
\end{lem}

\proof
Let us set
$\mu=
\pi_{0\sharp}(\eta)=
\pi_{1\sharp}(\eta)
$
and let
$$
X_1=\pi_0(D_1)\cap \pi_1(D_{1})\subset M
$$
be the set of points $x_0\in M$ such that a sequence
$x_{-1},x_0,x_1$ exists, with
$(x_{-1},x_0)\in D_{1}$
and $(x_0,x_1)\in D_{1}$.
Clearly, we have
$\mu(\pi_0(D_{1}))=\mu(\pi_1(D_{1}))=1$
hence $\mu(X_1)=1$.
Let 
$$
D_2= D_1\cap (X_1\times X_1)
$$
be the set of pairs $(x_0,x_1)\in M^2$ such that there exist
$x_{-1},x_0,x_1,x_2$ with $(x_i,x_{i+1})\in D_1$ for $i=-1,0,1$.
Let
$$
X_2=\pi_0(D_2)\cap \pi_1(D_{2})\subset M
$$
be the set of points $x_0\in M$ such that a sequence
$x_{-2},x_{-1},x_0,x_1,x_2$ exists, with
$(x_{i},x_{i+1})\in D_{1}$
for all $-2\leq i\leq 1$.
Since $\mu(X_1)=1$, we have $\eta(D_2)=1$,
hence $\mu(X_2)=1$.
By recurrence, we build a sequences $D_n\subset M\times M$ 
and $X_n\subset M$
such that 
$$
D_n=D_1\cap(X_{n-1}\times X_{n-1})
$$
and 
$$
X_n=\pi_0(D_n)\cap \pi_1(D_{n})\subset M.
$$
By recurrence, we see that $\eta(D_n)=1$
and that $\mu(X_n)=1$.
Now we have 
$$
D_{\infty}=\bigcap_{n\in \Zm} D_n
$$
hence $\eta(D_{\infty})=1$.
\qed

\subsection{Ergodic Construction}
It is worth explaining that the preceding
construction could have been performed in a quite different
way, which does not use our theory
of Kantorovich pairs, but relies on Ergodic theory,
as the first papers of Mather.

Consider $X=M^{\ZZ}$ endowed with the product  topology,
 so that $X$ is a metrizable compact space.
We shall denote by $\mM_T(X)$ the set of  Borelian probability measures on
$X$ which are invariant by translation.
More precisely, we denote 
 by $T:X\rightarrow
X$ the translation map
$$T(\ldots,a_{-2},a_{-1},a_0,a_1,a_{2},\ldots)
=(\ldots,b_{-2},b_{-1},b_0,b_1,b_{2},\ldots)$$
with $b_i=a_{i+1}$ for all $i\in \ZZ$,
so that $\mM_T(X)$ is the set of probability
measures $\nu$ on $X$ such that
$T_{\sharp}\nu=\nu$.

\begin{thm}\label{Mather2}
We have
$$
\min _{\nu\in \mM_T(X)}\int_X A(x_0,x_1)\, d\nu(x) =0.
$$
The measure $\nu$ is minimizing if and only 
if its marginal 
 $\eta=(\pi_0\times \pi_1)_{\sharp}\nu$
is concentrated on $D$.
\end{thm}

Note that Theorem \ref{Mather2} is equivalent to Theorem 
\ref{Mather} in view of the following:

\begin{lem} Let $\eta$ be a Borelian probability measure on $M^2$
such that 
$\pi_{0\sharp}(\eta)=\pi_{1\sharp}(\eta)$. Then there exists a Borelian
measure $\nu$ on $X$ that is $T$-invariant and such that
$\eta$ is its push-forward by the map $X\ni x\rightarrow (x_0,x_1)\in M^2$.
\end{lem}

\proof This follows from the Hahn-Kolmogorov extension theorem
(see e.g. Theorem 0.1.5 in \cite{Manebook}, Lemma 10.2.4
in \cite{Du:02} and Theorem 12.1.2 in \cite{Du:02}). 
Let $\Omega$ be the algebra
of finite unions of subsets $G$ of $X$ of the type $G=\Pi_{i\in \ZZ}G_i$ where
$G_i\neq M$ for at most a finite number of indices $i$ (the number
depending on $G$) and every $G_i$ is a Borelian subset of $M$. 
We first define the $T$-invariant probability measure $\nu$ on $\Omega$ 
and then apply
the Hahn-Kolmogorov extension theorem, which provides an unique extension to
the Borel  $\sigma$-algebra (by uniqueness, the extension is $T$-invariant).

Let $\eta=\int_{M}\eta_{x_1}d\mu(x_1)$ be the disintegration of $\eta$ with
respect to the projection $M^2\ni(x_0,x_1)\rightarrow x_1\in M$.
In particular $\mu=\pi_{1\sharp}(\eta)$ (see e.g.
Theorem 5.3.1 in \cite{AmGiSa} for the disintegration theorem).
Define for $m<n$
\begin{eqnarray*}
&&
\nu(\ldots\times M\times M\times G_{m}\times\ldots\times G_{n}
\times M\times M\times \ldots)\\
&=&\int_{G_m\times\ldots\times G_n}d\eta_{x_{m+1}}(x_{m})\ldots
d\eta_{x_{n}}(x_{n-1}) d\mu(x_n).
\end{eqnarray*}
This is well defined because if $G_{m-1}=M$ then
\begin{eqnarray*}
&&
\int_{G_{m-1}\times G_m\times\ldots\times G_n}d\eta_{x_{m}}(x_{m-1})
d\eta_{x_{m+1}}(x_{m})\ldots
d\eta_{x_{n}}(x_{n-1}) d\mu(x_n)\\
&=&\int_{G_m\times\ldots\times G_n}d\eta_{x_{m+1}}(x_{m})\ldots
d\eta_{x_{n}}(x_{n-1}) d\mu(x_n)
\end{eqnarray*}
and if $G_{n+1}=M$ then
\begin{eqnarray*}
&=&\int_{G_m\times\ldots\times G_n\times G_{n+1}}d\eta_{x_{m+1}}(x_{m})\ldots
d\eta_{x_{n}}(x_{n-1}) d\eta_{x_{n+1}}(x_{n}) d\mu(x_{n+1})\\
&=&\int_{G_n\times M}\left\{\int_{G_m\times\ldots\times G_{n-1}}
d\eta_{x_{m+1}}(x_{m})\ldots d\eta_{x_{n}}(x_{n-1})
\right\}d\eta_{x_{n+1}}(x_{n}) d\mu(x_{n+1})\\
&=&\int_{G_n\times M}\left\{\int_{G_m\times\ldots\times G_{n-1}}
d\eta_{x_{m+1}}(x_{m})\ldots d\eta_{x_{n}}(x_{n-1})
\right\}d\eta(x_{n},x_{n+1})\\
&=&\int_{G_n}\left\{\int_{G_m\times\ldots\times G_{n-1}}
d\eta_{x_{m+1}}(x_{m})\ldots d\eta_{x_{n}}(x_{n-1})
\right\} d\mu(x_{n})\\
&=&\int_{G_m\times\ldots\times G_n}d\eta_{x_{m+1}}(x_{m})\ldots
d\eta_{x_{n}}(x_{n-1}) d\mu(x_n)
\end{eqnarray*}
because $\mu=\pi_{1\sharp}(\eta)=\pi_{0\sharp}(\eta)$.

Clearly $\nu(X)=1$ and $\nu$ is $T$-invariant on $\Omega$.
\qed

Although we have proved the equivalence between 
Theorem \ref{Mather2} and
Theorem \ref{Mather}
we shall, as announced, detail another proof of Theorem \ref{Mather2}.

For $x\in X$ and every Borelian subset $B$, we define 
$$\tau_B(x)=\lim_{n\rightarrow +\infty} \frac 1 n \, \hbox{card}
\{0\leq j\leq n-1\,|\,T^j(x)\in B\}$$
(when the notation is used, it is understood that the limit exists). 
A Borelian probability $\nu$ on $X$ 
is ergodic if and only if,
for every Borelian subset $B\subset X$, there holds
$\tau_B(x)=\nu(B)$ $\nu$-almost surely.

Following Section II.6 in the book by Ma\~n\'e \cite{Manebook}, 
there exists
a Borel set $\Sigma\subset X$ such that $\nu(\Sigma)=1$
for each $\nu \in \mM_T(X)$,
and, for each $x\in \Sigma$,
the measure
$$
\nu_x:=
\lim_{n\rightarrow +\infty}\frac 1 n 
\sum_{j=0}^{n-1}\delta_{T^j(x)}
$$ 
is well defined and ergodic, where the limit is understood
in the sense of the weak topology, that is
\begin{equation}\label{eq: erg weak top} 
\forall f\in C(X,\RR)~~
\int_X f\,d\nu_x=\lim_{n\rightarrow \infty}\frac 1 n 
\sum_{j=0}^{n-1}f(T^j(x)).
\end{equation}
Moreover $\nu_x\in \mM_T(X)$ 
and $x$ belongs to the support of $\nu_x$ for all $x\in\Sigma$.
In addition, still following  \cite{Manebook}, we have
that the function $x\lmto \int fd\nu_x$ is $\nu$-integrable
and $T$-invariant,
and that
\begin{equation}\label{erg dec}
\int_X\left(\int_X fd\nu_x\right)d\nu=\int_X f\, d\nu.
\end{equation}
holds for
every $f\in \mL^1(X,\nu)$.
Note that the measure $\nu_x$ is the conditional probability measure
of $\nu$ with respect to the $\sigma$-algebra of $T$-invariant 
Borel sets.

We define the continuous function
 $\Gamma:X\rightarrow \RR$ by 
$\Gamma(x)= A(x_0,x_1)$.
By standard convexity arguments, the following minimum
is reached:
$$\alpha=\min_{\nu\in \mM_T(X)}\int_X \Gamma(x)d\nu(x).$$

Let us prove that $\alpha\leq 0$.
Fix $x_0\in M$. For all $\epsilon>0$, we can find
$n\geq 1$ and $x_1,\ldots,x_{n}\in M$
such that  
$$x_n=x_0\ \hbox{ and } \ \frac 1 n  \sum_{j=0}^{n-1}A(x_j,x_{j+1})<\epsilon$$
(thanks to assumption (A1)).
Let $x=(\ldots,x_0,\ldots,x_n,\ldots)\in X$ have periodic components 
with period $n$ and define $\nu\in \mM_T(X)$ by
$$\nu=\frac 1 n \sum_{j=0}^{n-1}\delta_{T^j(x)}$$
where  $\delta_{T^j(x)}$ is the Dirac measure at $T^j(x)$.
Then $\int_X\Gamma\, d\nu< \epsilon$, which proves that $\alpha\leq 0$ (because
$\epsilon$ can be chosen arbitrarily small).

Let $\nu\in \mM_T(X)$ be any optimal measure.
 The equality
$$
\int_{X}\left(\int_X \Gamma d\nu_x\right)d\nu=\int_X \Gamma\, d\nu=\alpha
$$
shows that 
$\int_X \Gamma d\nu_x=\alpha$ for $\nu$-almost all $x\in \Sigma$.
For such a $x$, we get
\begin{equation}\label{eq: 0 limit}
0\geq \alpha=\lim_{n\rightarrow +\infty}\frac 1 n \sum_{j=0}^{n-1}\Gamma(T^j(x))
=\lim_{n\rightarrow +\infty}\frac 1 n \sum_{j=0}^{n-1}A(x_j,x_{j+1}).
\end{equation}
Assume for a while that $x_0\not\in \mA$. Then there exists a neighborhood
$U$ of $x_0$ in $M$, $\delta>0$ and $N\geq 1$ such that
\begin{equation}\label{eq: def of N}
\hbox{$A_n(y_0,z_0)>\delta>0$ for all $y_0,z_0\in U$ and $n\geq N$}
\end{equation}
(we use here the equicontinuity of the functions $A_n$).
Setting $\widetilde U=\{y\in X\,|\,y_0\in U\}$, we get
$$
0<\nu_x(\widetilde U)
\leq \liminf_{n\rightarrow +\infty} \frac 1 n \, \hbox{card}
\{0\leq j\leq n-1\,|\,x_j\in U\}.
$$
The first inequality is a consequence of the fact that $x$ is in the support
of $\nu_x$ and the second one follows from \eqref{eq: erg weak top} 
and the fact that
the characteristic function of $U$ is the supremum of an increasing sequence
of continuous functions.
We denote by  $(x_{j_k}:k\geq 0)$ the sequence of components 
of $x$ in $U$ (of non-negative index).
We obtain (see \eqref{eq: def of N})
$$
0<\nu_x(\widetilde U)
\leq \liminf_{m\rightarrow +\infty} \frac {mN} {j_{mN}} 
$$
and the contradiction
\begin{eqnarray*}
&&\liminf_{m\rightarrow \infty}
\frac 1 {j_{mN}} \sum_{j=0}^{j_{mN}-1}A(x_j,x_{j+1})
= \liminf_{m\rightarrow \infty}
\frac 1 {j_{mN}}\sum_{k=0}^{m-1} \sum_{j=j_{kN}}^{j_{(k+1)N}-1}A(x_j,x_{j+1})
\\&&\geq \liminf_{m\rightarrow \infty}
\frac 1 {j_{mN}}\sum_{k=0}^{m-1} A_{j_{(k+1)N}-j_{kN}}
(x_{j_{kN}},x_{j_{(k+1)N}})
\\&& \geq 
\liminf_{m\rightarrow \infty}
\frac {m\delta} {j_{mN}}\geq \nu_x(\widetilde U)\delta/N\, >\, 0
\end{eqnarray*}
(compare with \eqref{eq: 0 limit}).
This contradiction shows that $x_0\in \mA$ for $\nu$-almost all $x$,
that is, the marginal 
$\mu=\pi_{0\sharp} \nu$
is concentrated on $\mA$.

Let us now check that $\alpha\geq 0$. For contradiction, suppose
$\alpha<0$. Then
for $x\in \Sigma$ as above
such that $\nu_x\in \mM_T(X)$ and $\Gamma(\nu_x)=\alpha$, we get 
\begin{eqnarray*}
0&>&\alpha=\Gamma(\nu_x)
=\lim_{n\rightarrow +\infty}\frac 1 n \sum_{j=0}^{n-1}A(x_j,x_{j+1})
\\&=&\lim_{n\rightarrow +\infty}\frac 1 {n+1} 
\left(\sum_{j=0}^{n-1}A(x_j,x_{j+1})+A(x_n,x_0)\right)
\\& \geq& \limsup_{n\rightarrow +\infty}\frac 1 {n+1} A_{n+1}(x_0,x_{0}).
\end{eqnarray*}
This contradicts $l=0$ (see hypothesis (A1)).

We have proved that $\alpha=0$, and that every minimizing $T$-invariant measure
$\nu$ has its  marginal $\mu=\pi_{0\sharp} \nu$
concentrated on $\mA$.
Let us now  prove that every minimizing measure
$\nu\in \mM_T(X)$  is supported on 
$\{x\in X\,|\,A(x_0,x_1)+c(x_1,x_0)=0\}.$
Let $x$ belong to the support of $\nu$ and observe that
(see Lemma \ref{propdec})
$c(x_1,y_1)\leq c(x_1,y_0)+A(y_0,y_1)$ for all $y_0,y_1\in M$. Therefore
$$0=\int_X c(x_1,y_1)-c(x_1,y_0)   \, d\nu(y)
\leq \int_X A(y_0,y_1)   \, d\nu(y)=\alpha=0$$
and
$$\int_X A(y_0,y_1)- c(x_1,y_1)+c(x_1,y_0)   \, d\nu(y)=0$$
where the integrand is non negative.
Hence $c(x_1,y_1)= c(x_1,y_0)+A(y_0,y_1)$ for $\nu$-almost all $y$.
Since $x$ is in the support of $\nu$, we get
$c(x_1,x_1)= c(x_1,x_0)+A(x_0,x_1)$.
We have just seen that $y_0\in \mA$ for $\nu$-almost all $y$.
By the  $T$-invariance of $\nu$, we also have
$y_1\in \mA$ for $\nu$-almost all $y$. Since  $x$ is in the support
of $\nu$, we therefore obtain $x_1\in \mA$ and
$0=c(x_1,x_1)= c(x_1,x_0)+A(x_0,x_1)$.

Finally let $\nu\in \mM_T(X)$ be concentrated on
$$\widetilde D=\{y\in X\,|\,A(y_0,y_1)+c(y_1,y_0)=0\}$$
and let us prove that $\int_X \Gamma d\nu=\alpha$.
By \eqref{erg dec} applied to the characteristic function of $\widetilde D$,
we get that $\nu_x(\widetilde D)=1$ for $\nu$-almost all $x\in \Sigma$.
By \eqref{erg dec} applied to $\Gamma$, we see that it suffices
to check that $\int_X \Gamma d\nu_x=\alpha$ for all $x\in \Sigma$
such that $\nu_x$ is concentrated on $\widetilde D$.
This follows from (C1):
\begin{multline*}
0=\alpha\leq \int_X A(y_0,y_1)\, d\nu_x(y)
=-\int_X c(y_1,y_0)\, d\nu_x(y)
\\=-\lim_{n\rightarrow \infty}\frac 1 n \sum_{j=0}^{n-1}c(x_{j+1},x_{j})
\leq -\liminf_{n\rightarrow \infty}\frac 1 n c(x_n,x_0)
=0.
\end{multline*}
\qed

\section{Aubry-Mather theory}
We now briefly explain the  relations between our
discussions and the literature on Aubry-Mather theory,
and especially \cite{Ma:91}, \cite{Mane} and \cite{Fa:un}.
From now on, the space $M$ is a compact connected manifold and
we consider a $C^2$ Lagrangian function 
$L:TM\times \Rm\lto \Rm$
as in the Introduction.
 In this context, we define $A:M\times M\lto \Rm$
 by 
 $$
 A(x,y)=\min_{\gamma}
 \int_0^1 L(\gamma(t),\dot \gamma(t),t)dt,
 $$
 where the minimum is taken on the set of $C^2$ curves
 $\gamma:[0,1]\lto M$ which satisfy 
 $\gamma(0)=x$ and $\gamma(1)=y$.
 
 The function $c$ defined by (\ref{defc})
is one of the central objects of Mather's theory
of globally minimizing orbits, see \cite{Ma:91}.
He called it the Peierls barrier.
It contains most of the information concerning
the globally minimizing orbits, as was explained by Mather,
see also \cite{symp}.
The set $\mA$ of points $x\in M$ such that $c(x,x)=0$
is called the projected Aubry set. 
It is especially important because Mather
proved the existence of a vectorfield $X(x)$
on $\mA$ whose graph is invariant under the 
Lagrangian flow $\varphi$.
This invariant set is called the Aubry set.
The analog of the Aubry set in our general theory
is the set $D$ defined in the beginning of section 
\ref{dyn}.

The operators $T^\pm$
have  been introduced by Albert Fathi in this context,
see \cite{Fa:97a},\cite{Fa:97b} and  \cite{Fa:un}.
He called Weak KAM solutions the fixed points of $T^-$,
and we call backward weak KAM solutions the fixed points of $T^+$.
He also noticed that, for each weak KAM solution
$\phi_1$, there exists one and only one backward weak KAM
solution $\phi_0$ which is equal to $\phi_1$
on the projected Aubry set. This is the main part
of our Theorem \ref{pairs}.
Albert Fathi also proved Lemma %
\ref{Kanto}
in this context. 
Our novelty in these matters consists of pointing out
and using the equivalence with Kantorovich admissible
pairs, which allows, for example, a strikingly simple proof of the important 
result of Fathi called 
Lemma %
\ref{Kanto}
in our paper.
The representation of weak KAM solutions given in 
Proposition \ref{admissiblepairs} was obtained by
Contreras in \cite{Co:01}.

The minimizing measures of Theorem \ref{Mather}
are the famous Mather  measures, see \cite{Ma:91}.
To be more precise, we should say that there is a natural 
bijection between the set of minimizing measures in 
Theorem \ref{Mather} and the set of Mather measures.
This bijection is described in \cite{BeBu:first}.
In order to give the reader a clue of this bijection,
let us recall that the Mather measures are probability measures
on the tangent bundle $TM$, and that the minimizing measures
of Theorem \ref{Mather} are probability measures on $M\times M$.
Denoting by $\varphi$ the time-one Lagrangian flow,
and by $\pi:TM\lto M$ the standard projection,
we have a well-defined mapping 
$(\pi,\pi\circ \varphi)_{\sharp}$ from the set of probability
measures on $TM$ to the set of probability measures  
on $M\times M$.
This mapping induces a bijection between the set of Mather
measures on $TM$ and the set of 
minimizing measures of Theorem \ref{Mather}.

The part of Theorem \ref{Mather} stating that the minimizing
measures are precisely the measures supported on $D$
is the analogous in our setting of the theorem of Ma\~n\'e
stating that all invariant measures supported on the Aubry set 
are minimizing, see \cite{Ma:92}.

\end{document}